\documentclass[10pt,%
]{bmc_article}

\makeatletter
\renewcommand\maketitle{%
   \begin{flushleft}\mbox{%
    \global\let\@date\@empty% @date = null.
    {\sffamily\begin{minipage}{\textwidth}%
             \@maketitle%
             {\raggedright%
                {\noindent\@address}\\ \hbox{}%
                {\noindent\@bmc@email}\\ \hbox{}%
                {\noindent\@corresponding}%
             }%end \raggedright
      \end{minipage}%
      }%
    \renewcommand\thefootnote{\old@thefootnote}%
    }% \mbox
    \end{flushleft}%
}
\makeatother
\usepackage{cite} % Make references as [1-4], not [1,2,3,4]
\usepackage{url}  % Formatting web addresses
\usepackage{ifthen}  % Conditional
\usepackage{multicol}   %Columns
\usepackage{amsmath}
\usepackage{mathrsfs}
\usepackage{latexsym}
\usepackage{algorithm} % remove when submit
\usepackage{algorithmic}
\usepackage{multirow}
\usepackage{color}
\newcommand{\STATECOMMENT}[1]{\STATE$\rhd$ #1}
\newcommand{\MYPRINT}{\STATE\textbf{print}\ \ }
\newcommand{\MYRETURN}{\STATE\textbf{return}\ \ }
\newcommand{\GOTO}{\STATE\textbf{goto}\ \ }

\DeclareMathOperator{\foldpk}{\texttt{Cross}}
\DeclareMathOperator{\checkstru}{\textsc{Check-Stru}}
\DeclareMathOperator{\makestart}{\textsc{Make-Start}}
\DeclareMathOperator{\adjustseq}{\textsc{Adjust-Seq}}
\DeclareMathOperator{\decompose}{\textsc{Decompose}}
\DeclareMathOperator{\aw}{\textsc{Local-Search}}

\DeclareMathOperator{\bigo}{O} \DeclareMathOperator{\crosscomp}{c}

\def\baseA{{\bf A}}
\def\baseC{{\bf C}}
\def\baseG{{\bf G}}
\def\baseU{{\bf U}}
\def\pairAU{\textbf{A-U}}
\def\pairCG{\textbf{C-G}}
\def\pairGC{\textbf{G-C}}
\def\pairGU{\textbf{G-U}}
\def\pairUA{\textbf{U-A}}
\def\pairUG{\textbf{U-G}}

\newtheorem{theorem}{Theorem}

\usepackage{graphicx} % added by James, remove when submit
\DeclareGraphicsExtensions{.eps} % added by James
\graphicspath{{fig/}} % added by James

\newboolean{publ}

\newenvironment{bmcformat}{%\begin{raggedright}%Linda changed
\baselineskip20pt\sloppy\setboolean{publ}{false}}{%\end{raggedright}
\baselineskip20pt\sloppy}
\begin{document}
\begin{bmcformat}
%%%
%%%%%%%%%%%%%%%%%%%%%%%%%%%%%%%%%%%%%%%%%%%%%%%%%%%%%%%%%%%%%%%%%%%%%%%%%%%%
%%%
\title{Inverse Folding of RNA Pseudoknot Structures}

\author{%
  James Z.M. Gao$^1$ \email{Gao: gzm55@cfc.nankai.edu.cn}%
\and
  Linda Y.M. Li$^1$ \email{Li: liyanmei@mail.nankai.edu.cn}%
and
  Christian M. Reidys$^1$\correspondingauthor%
    \email{\correspondingauthor Reidys: duck@santafe.edu}%
}

\address{%
  \iid (1) Center for Combinatorics, LPMC-TJKLC, Nankai University,
           Tianjin 300071, PR China
}

\maketitle
%%%
%%%%%%%%%%%%%%%%%%%%%%%%%%%%%%%%%%%%%%%%%%%%%%%%%%%%%%%%%%%%%%%%%%%%%%%%%%%%%
%%%

%%%
%%%%%%%%%%%%%%%%%%%%%%%%%%%%%%%%%%%%%%%%%%%%%%%%%%%%%%%%%%%%%%%%%%%%%%%%%%%%
%%%
\begin{abstract}
  \paragraph{Background:}
  RNA exhibits a variety of structural configurations.
Here we consider a structure
  to be tantamount to the noncrossing Watson-Crick and \pairGU-base pairings
  (secondary structure) and additional cross-serial base pairs. These interactions
  are called pseudoknots and are observed across the whole spectrum of RNA
  functionalities. In the
  context of studying natural RNA structures, searching for new ribozymes and
  designing artificial RNA, it is of interest to find RNA sequences folding
  into a specific structure and to analyze their induced neutral networks.
  Since the established inverse folding algorithms, {\tt RNAinverse}, {\tt
  RNA-SSD} as well as {\tt INFO-RNA} are limited to RNA secondary structures,
  we present in this paper the inverse folding algorithm {\tt Inv} which can
  deal with $3$-noncrossing, canonical pseudoknot structures.

%%%%%%%%%%%%%%%%%%%%%%%%%%%%%%%%%%%%%%%%%%%%%%%%%%%%%%%%%%%%%%%%%%%%%%%%%%%%%%%%
  \paragraph{Results:}

  In this paper we present the inverse folding algorithm {\tt Inv}. We give a
  detailed analysis of {\tt Inv}, including pseudocodes.
  We show that {\tt Inv} allows to design in particular $3$-noncrossing nonplanar
  RNA pseudoknot $3$-noncrossing RNA structures--a class which is difficult to
  construct via dynamic programming routines.
  {\tt Inv} is freely available at
  \url{http://www.combinatorics.cn/cbpc/inv.html}.

%%%%%%%%%%%%%%%%%%%%%%%%%%%%%%%%%%%%%%%%%%%%%%%%%%%%%%%%%%%%%%%%%%%%%%%%%%%%%%%%
  \paragraph{Conclusions:}

  The algorithm {\tt Inv} extends inverse folding capabilities to RNA
  pseudoknot structures. In comparison with {\tt RNAinverse} it uses new ideas,
  for instance by considering sets of competing structures.
  As a result, {\tt Inv}
  is not only able to find novel sequences even for RNA secondary structures, it
  does so in the context of competing structures that potentially exhibit
  cross-serial interactions.

\end{abstract}

%%%
%%%%%%%%%%%%%%%%%%%%%%%%%%%%%%%%%%%%%%%%%%%%%%%%%%%%%%%%%%%%%%%%%%%%%%%%%%%%%%%%
%%%

\ifthenelse{\boolean{publ}}{\begin{multicols}{2}}{}

\twocolumn % remove when submit

%
%%%%%%%%%%%%%%%%%%%%%%%%%%%%%%%%%%%%%%%%%%%%%%%%%%%%%%%%%%%%%%%%%%%%%%%%%%%%%%%%
%%%
\section{Introduction}
%%%
%%%%%%%%%%%%%%%%%%%%%%%%%%%%%%%%%%%%%%%%%%%%%%%%%%%%%%%%%%%%%%%%%%%%%%%%%%%%%%%%
%%%

Pseudoknots are structural elements of central importance in RNA structures
\cite{Westhof:92}, see Figure~\ref{F:glmS}. They represent cross-serial base
pairing interactions between RNA nucleotides that are functionally important in
tRNAs, RNaseP \cite{Loria:96}, telomerase RNA \cite{Butcher}, and ribosomal
RNAs \cite{Konings:95a}. Pseudoknot structures are being observed in the
mimicry of tRNA structures in plant virus RNAs as well as the binding to the
HIV-1 reverse transcriptase in {\it in vitro} selection experiments
\cite{Tuerk:92}. Furthermore basic mechanisms, like ribosomal frame shifting,
involve pseudoknots \cite{Chamorro:91a}.

%%%
%%%%%%%%%%%%%%%%%%%%%%%%%%%%%%%%%%%%%%%%%%%%%%%%%%%%%%%%%%%%%%%%%%%%%%
%%%
\begin{figure}
  \centering
  \includegraphics[width=\columnwidth]{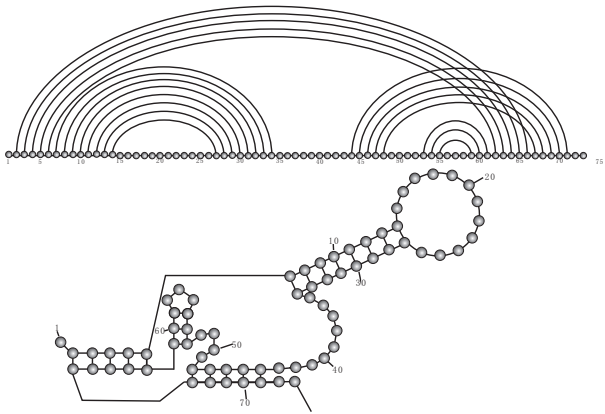}
  \caption{\small
    The pseudoknot structure of the glmS ribozyme pseudoknot P1.1
     \cite{W:glms}  as a diagram (top) and as a planar graph (bottom).
   }
  \label{F:glmS}
\end{figure}
%%%
%%%%%%%%%%%%%%%%%%%%%%%%%%%%%%%%%%%%%%%%%%%%%%%%%%%%%%%%%%%%%%%%%%%%%%%%%%%%%%%%
%%%
%%%
%%%
%%%%%%%%%%%%%%%%%%%%%%%%%%%%%%%%%%%%%%%%%%%%%%%%%%%%%%%%%%%%%%%%%%%%%%%
%%%

Despite them playing a key role in a variety of contexts, pseudoknots are
excluded from large-scale computational studies. Although the problem has
attracted considerable attention in the last decade, pseudoknots are considered
a somewhat ``exotic'' structural concept. For all we know \cite{Lyngso}, the
{\it ab initio} prediction of general RNA pseudoknot structures is NP-complete
and algorithmic difficulties of pseudoknot folding are confounded by the fact
that the thermodynamics of pseudoknots is far from being well understood.

As for the folding of RNA secondary structures, Waterman {\it et al}
\cite{Waterman:78a, Waterman:86}, Zuker {\it et al} \cite{Zuker:1981}
and Nussinov \cite{Nussinov:1980}
established the dynamic programming ({\tt DP}) folding routines. The
first mfe-folding algorithm for RNA secondary structures, however,
dates back to the 60's \cite{Jaces:1960, Tinoco:1971, Delisi:1971}.
For restricted classes of pseudoknots, several algorithms have been
designed: Rivas and Eddy \cite{RE:98}, Dirks and
Pierce\cite{Dirks:2004}, Reeder and Giegerich \cite{Reeder:04} and
Ren {\it et al} \cite{Ren:05}.
Recently, a novel {\it ab initio}
folding algorithm {\tt Cross} has been introduced
\cite{Reidys:08algo}. ${\tt Cross}$ generates minimum free energy
(mfe), $3$-noncrossing, $3$-canonical RNA structures,
i.e.~structures
that do not contain three or more mutually crossing arcs and in
which each stack, i.e.~sequence of parallel arcs, see eq.~(\ref{E:cano}),
has size greater or equal than three.
In particular, in a $3$-canonical structure there
are no isolated arcs, see Figure~\ref{F:canonical}.

%%%
%%%%%%%%%%%%%%%%%%%%%%%%%%%%%%%%%%%%%%%%%%%%%%%%%%%%%%%%%%%%%%%%%%%%%%%
%%%
\begin{figure}[ht]
\centerline{\includegraphics[width=0.4\textwidth]{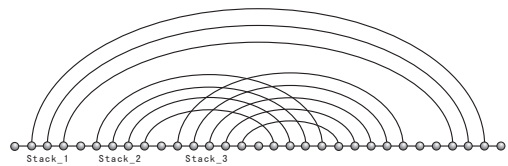}}
\caption{\small $\sigma$-canonical
RNA structures: each stack of
``parallel'' arcs has to have minimum size $\sigma$. Here we display a
$3$-canonical structure.}
\label{F:canonical}
\end{figure}
%%%
%%%%%%%%%%%%%%%%%%%%%%%%%%%%%%%%%%%%%%%%%%%%%%%%%%%%%%%%%%%%%%%%%%%%%%
%%%
%%%

The notion of mfe-structure is based on a
specific concept of pseudoknot loops and respective loop-based
energy parameters. This thermodynamic model was conceived by Tinoco
and refined by Freier, Turner, Ninio, and others \cite{Tinoco:1971,
Borer:1974, Papanicolaou:1984, Turner:1988, Walter:1994, Xia:1998}.

%%%
%%%%%%%%%%%%%%%%%%%%%%%%%%%%%%%%%%%%%%%%%%%%%%%%%%%%%%%%%%%%%%%%%%%%%%%%%%%%%%%%
%%%
\subsection{$k$-noncrossing, $\sigma$-canonical RNA pseudoknot structures}
%%%
%%%%%%%%%%%%%%%%%%%%%%%%%%%%%%%%%%%%%%%%%%%%%%%%%%%%%%%%%%%%%%%%%%%%%%%%%%%%%%%%
%%%

Let us turn back the clock: three decades ago Waterman {\it et al.}
\cite{Waterman:79a}, Nussinov {\it et al.} \cite{Nussinov:1980}
and Kleitman {\it et al.} in \cite{Kleitman:70} analyzed RNA secondary
structures. Secondary structures are coarse grained RNA contact structures, see
Figure~\ref{F:tRNA1}.

%%%
%%%%%%%%%%%%%%%%%%%%%%%%%%%%%%%%%%%%%%%%%%%%%%%%%%%%%%%%%%%%%%%%%%%%%%%
%%%
\begin{figure}[ht]
\centerline{\includegraphics[width=0.5\textwidth]{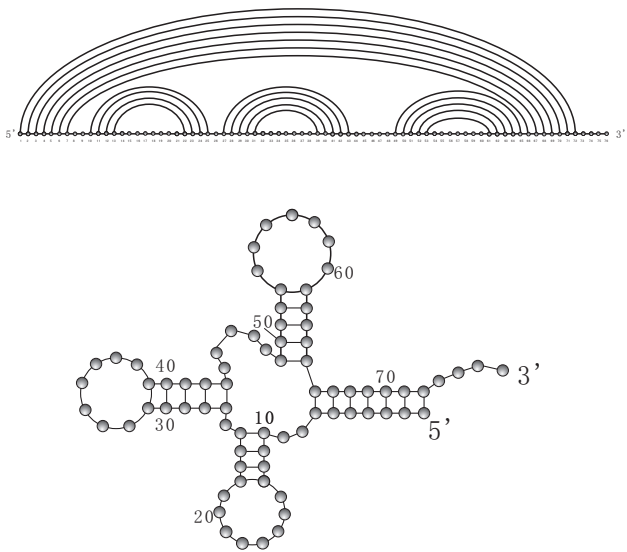}}
\caption{\small The phenylalanine tRNA secondary structure represented
as $2$-noncrossing diagram (top) and as planar graph (bottom).}
\label{F:tRNA1}
\end{figure}
%%%
%%%%%%%%%%%%%%%%%%%%%%%%%%%%%%%%%%%%%%%%%%%%%%%%%%%%%%%%%%%%%%%%%%%%%%
%%%
Secondary structures can be represented as diagrams, i.e.~labeled graphs over
the vertex set $[n]=\{1, \dots, n\}$ with vertex degrees $\le 1$, represented
by drawing its vertices on a horizontal line and its arcs $(i,j)$ ($i<j$), in
the upper half-plane, see Figure~\ref{F:glmS}~and Figure~\ref{F:secondary}.

Here, vertices and arcs correspond to the nucleotides \baseA, \baseG, \baseU,
\baseC \ and
Watson-Crick (\pairAU, \pairGC) and (\pairUG) base pairs, respectively.

In a diagram, two arcs $(i_1,j_1)$ and $(i_2,j_2)$ are called crossing if
$i_1<i_2<j_1<j_2$ holds. Accordingly, a $k$-crossing is a sequence of arcs
$(i_1,j_1),\dots,(i_k,j_k)$ such that $i_1<i_2<\dots<i_k<j_1<j_2<\dots <j_k$,
see Figure~\ref{F:p_cross}.

%%%
%%%%%%%%%%%%%%%%%%%%%%%%%%%%%%%%%%%%%%%%%%%%%%%%%%%%%%%%%%%%%%%%%%%%%%%
%%%
\begin{figure}[ht]
\centerline{\includegraphics[width=0.4\textwidth]{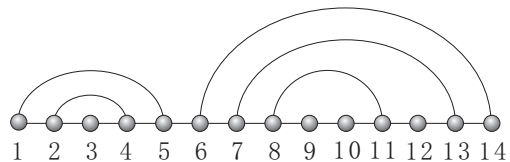}}
\caption{\small 
Setting $k=2$ we observe that secondary structures are
a particular type of $k$-noncrossing structures. They coincide with
noncrossing diagrams having minimum arc-length two.}
\label{F:secondary}
\end{figure}
%%%
%%%%%%%%%%%%%%%%%%%%%%%%%%%%%%%%%%%%%%%%%%%%%%%%%%%%%%%%%%%%%%%%%%%%%%
%%%
%%%
%%%%%%%%%%%%%%%%%%%%%%%%%%%%%%%%%%%%%%%%%%%%%%%%%%%%%%%%%%%%%%%%%%%%%%%
%%%
\begin{figure}[ht]
\centerline{\includegraphics[width=0.4\textwidth]{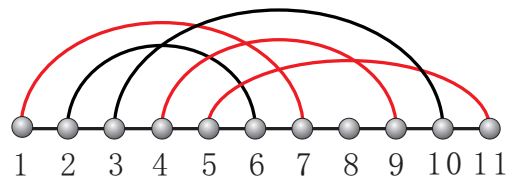}}
\caption{\small
$k$-noncrossing diagrams: we display
a $4$-noncrossing diagram containing the three mutually
crossing arcs $(1,7),(4,9),(5,11)$ (drawn in red).}
\label{F:p_cross}
\end{figure}
%%%
%%%%%%%%%%%%%%%%%%%%%%%%%%%%%%%%%%%%%%%%%%%%%%%%%%%%%%%%%%%%%%%%%%%%%%
%%%
We call diagrams containing at most $(k-1)$-crossings, $k$-noncrossing
diagrams. RNA secondary structures have no crossings in their diagram
representation, see Figure~\ref{F:tRNA1} and Figure~\ref{F:secondary}, and
are therefore $2$-noncrossing diagrams.
A structure in which any stack has at least
size $\sigma$ is called $\sigma$-canonical, where a stack\index{stack} of
size $\sigma$ is a sequence of ``parallel'' arcs of the form
\begin{equation}\label{E:cano}
((i,j),(i+1,j-1),\dots,(i+(\sigma-1),j-(\sigma-1))).
\end{equation}

As a natural generalization of RNA secondary structures $k$-noncrossing RNA
structures \cite{Reidys:07pseu,Reidys:08lego,Reidys:pnas09} were introduced. A
$k$-noncrossing RNA structure is $k$-noncrossing diagram without arcs of the
form $(i,i+1)$.  In the following we assume $k=3$, i.e.~in the diagram
representation there are at most two mutually crossing arcs, a minimum
arc-length of four and a minimum stack-size of three base pairs.  The notion
$k$-noncrossing stipulates that the complexity of a pseudoknot is related to
the maximal number of mutually crossing bonds. Indeed, most natural RNA
pseudoknots are $3$-noncrossing \cite{Stadler:99}.

%%%
%%%%%%%%%%%%%%%%%%%%%%%%%%%%%%%%%%%%%%%%%%%%%%%%%%%%%%%%%%%%%%%%%%%%%%%%%%%%%%%
%%%
\subsection{Neutral networks}\label{S:neutral}
%%%
%%%%%%%%%%%%%%%%%%%%%%%%%%%%%%%%%%%%%%%%%%%%%%%%%%%%%%%%%%%%%%%%%%%%%%%%%%%%%%%
%%%

Before considering an inverse folding algorithm into specific RNA structures
one has to have at least some rationale as to why there exists
{\it one} sequence realizing a given target as mfe-configuration.
In fact this is, on the level of
entire folding maps, guaranteed by the combinatorics of the target structures
alone.
It has been shown in \cite{Reidys:08ma}, that the numbers of $3$-noncrossing RNA
pseudoknot structures, satisfying the biophysical constraints
grows asymptotically as
$c_3 n^{-5} 2.03^n$, where $c_3>0$ is some explicitly known constant.
In view of the central limit theorems of \cite{Fenix:08}, this fact implies
the existence of extended (exponentially large) sets of sequences that all
fold into one $3$-noncrossing RNA pseudoknot structure, $S$.
In other words, the combinatorics of $3$-noncrossing RNA structures alone
implies that there are many sequences mapping (folding) into a single
structure. The set of all such sequences is called the neutral
network\footnote{the term ``neutral network'' as opposed to ``neutral set''
stems from giant component results of random induced subgraphs of $n$-cubes.
That is, neutral networks are typically connected in sequence space} of
the structure $S$ \cite{Reidys:97generic,Reidys:08local}, see
Figure~\ref{F:network}.

%%%
%%%%%%%%%%%%%%%%%%%%%%%%%%%%%%%%%%%%%%%%%%%%%%%%%%%%%%%%%%%%%%%%%%%%%%%%%%%%%%%
%%%
\begin{figure}
  \centering
  \includegraphics[width=\columnwidth]{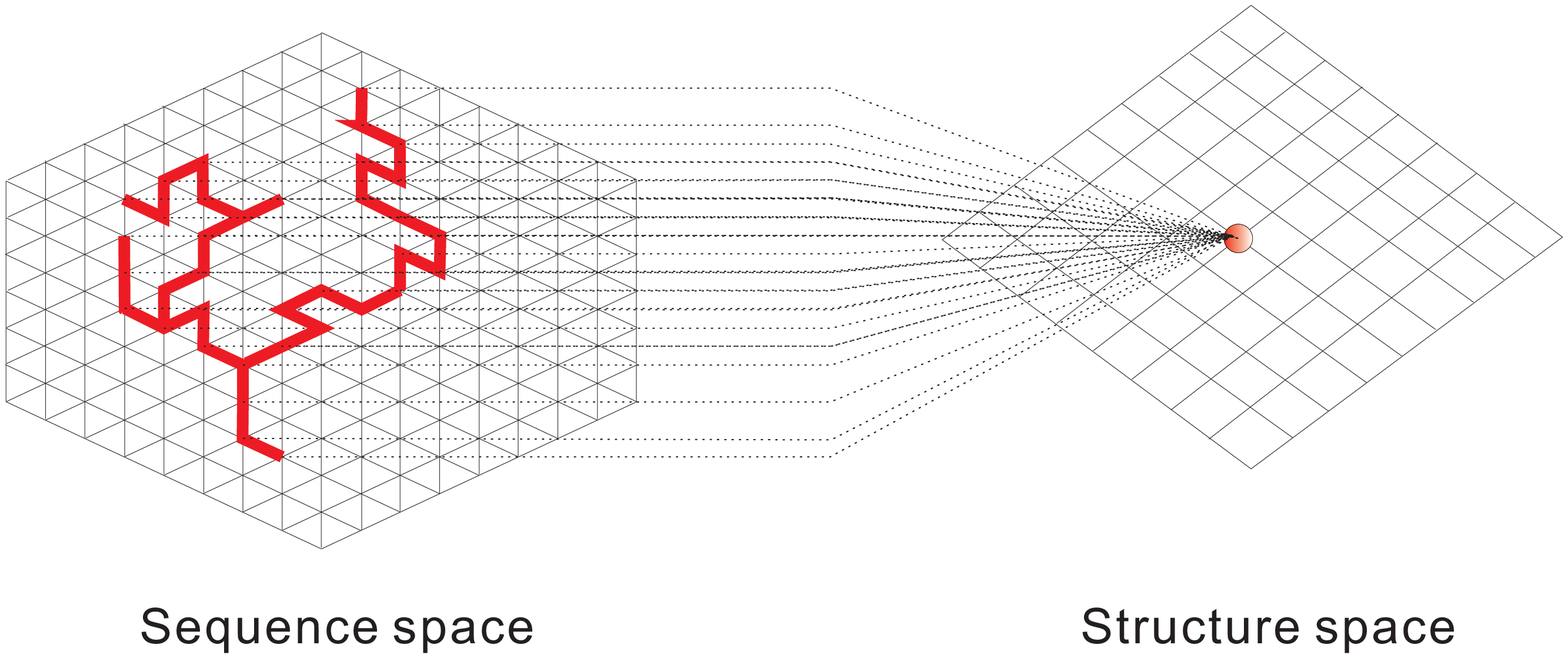}
  \caption{\small
  Neutral networks in sequence space: we display sequence space (left) and
  structure space (right) as grids. We depict a set of sequences that all fold into
  a particular structure. Any two of these sequences are connected by a red edge.
  The neutral network of this fixed structure consists of all sequences folding
  into it and is typically a connected subgraph of sequence space. }
  \label{F:network}
\end{figure}
%%%
%%%%%%%%%%%%%%%%%%%%%%%%%%%%%%%%%%%%%%%%%%%%%%%%%%%%%%%%%%%%%%%%%%%%%%%%%%%%%%%
%%%

By construction, all the sequences contained in such a neutral network are all
compatible with $S$. That is, at any two positions paired in $S$, we
find two bases capable of forming a bond (\pairAU, \pairUA, \pairGC,
\pairCG, \pairGU{} and \pairUG), see Figure~\ref{F:compatible}. Let
$s'$ be a sequence derived via a mutation\footnote{note: we do not
consider insertions or deletions.} of $s$. If $s'$ is again
compatible with $S$, we call this mutation ``compatible''.

%%%
%%%%%%%%%%%%%%%%%%%%%%%%%%%%%%%%%%%%%%%%%%%%%%%%%%%%%%%%%%%%%%%%%%%%%%%%%%%%%%%
%%%
\begin{figure}
  \centering
  \includegraphics[width=\columnwidth]{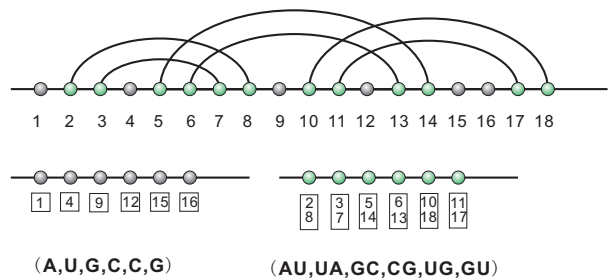}
  \caption{\small  A structure
  and a particular compatible sequence organized in
  the segments of unpaired and paired bases.}
  \label{F:compatible}
\end{figure}
%%%
%%%%%%%%%%%%%%%%%%%%%%%%%%%%%%%%%%%%%%%%%%%%%%%%%%%%%%%%%%%%%%%%%%%%%%%%%%%%%%%
%%%
Let $C[S]$ denote the set of $S$-compatible sequences. The structure
$S$ motivates to consider a new adjacency relation
within $C[S]$.
Indeed, we may
reorganize a sequence $(s_1,\dots,s_n)$ into the pair
\begin{equation}\label{E:decompose}
  \left((u_1,\dots,u_{n_u}),(p_1,\dots,p_{n_p})\right),
\end{equation}
where the $u_h$ denotes the unpaired nucleotides and the $p_h=(s_{i},s_{j})$
denotes base pairs, respectively, see Figure~\ref{F:compatible}. We can then
view $s_{u}=(u_1,\dots,u_{n_u})$ and $s_{p}=(p_1,\dots,p_{n_p})$ as elements of
the formal cubes $Q_4^{n_u}$ and $Q_6^{n_p}$, implying the new adjacency
relation for elements of $C[S]$.

Accordingly, there are two types of compatible neighbors in the sequence space
${\sf u}$- and ${\sf p}$-neighbors: a ${\sf u}$-neighbor has Hamming distance
one and differs exactly by a point mutation at an unpaired position.
Analogously a ${\sf p}$-neighbor differs by a compensatory base pair-mutation,
see Figure~\ref{F:mutate}.

%%%
%%%%%%%%%%%%%%%%%%%%%%%%%%%%%%%%%%%%%%%%%%%%%%%%%%%%%%%%%%%%%%%%%%%%%%%%%%%%%%%
%%%
\begin{figure}
  \centering
  \includegraphics[width=\columnwidth]{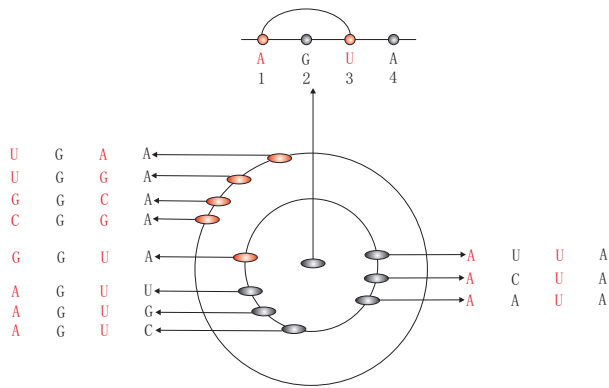}
  \caption{\small
  Diagram representation of an RNA structure (top) and
  its induced compatible neighbors in sequence space (bottom).
  Here the neighbors on the inner circle have Hamming distance one
  while those on the outer circle have Hamming distance two. Note that
  each base pair gives rise to five compatible neighbors (red) exactly one of
  which being in Hamming distance one. }
  \label{F:mutate}
\end{figure}
%%%
%%%%%%%%%%%%%%%%%%%%%%%%%%%%%%%%%%%%%%%%%%%%%%%%%%%%%%%%%%%%%%%%%%%%%%%%%%%%%%%
%%%

Note, however, that a {\sf p}-neighbor has either Hamming distance one
($\pairGC \mapsto \pairGU$) or Hamming distance two ($\pairGC \mapsto
\pairCG$). We call a {\sf u}- or a {\sf p}-neighbor, $y$, a compatible
neighbor. In light of the adjacency notion for the set of compatible
sequences we call the set of all sequences folding into $S$ the
neutral network of $S$.
By construction, the neutral network of $S$ is
contained in $C[S]$.
If $y$ is contained in the neutral network we
refer to $y$ as a neutral neighbor. This gives rise to consider the
compatible and neutral distance of the two sequences,
denoted by $C(s,s')$ and $N(s,s')$.
These are the minimum length of a $C[S]$-path and path in the neutral
network between $s$ and $s'$, respectively.
Note that since each neutral
path is in particular a compatible path, the compatible distance is
always smaller or equal than the neutral distance.

In this paper we study the inverse folding problem for RNA pseudoknot
structures: for
a given $3$-noncrossing target structure $S$,
we search for sequences from $C[S]$, that have $S$ as mfe configuration.

%%%
%%%%%%%%%%%%%%%%%%%%%%%%%%%%%%%%%%%%%%%%%%%%%%%%%%%%%%%%%%%%%%%%%%%%%%%%%%%%%%%%
%%%
\section{Background}\label{S:back}
%%%
%%%%%%%%%%%%%%%%%%%%%%%%%%%%%%%%%%%%%%%%%%%%%%%%%%%%%%%%%%%%%%%%%%%%%%%%%%%%%%%%
%%%

For RNA secondary structures, there are three different strategies for inverse
folding, {\tt RNAinverse},
{\tt RNA-SSD} and {\tt INFO-RNA} \cite{Hofacker:1994,SSD:2004,Busch:algo},

They all generate via a local search routine iteratively sequences, whose
structures have smaller and smaller distances to a given target.
Here the distance between two structures
is obtained by aligning them as diagrams and counting
``$0$'', if a given position is either unpaired or incident to an arc
contained in both structures and ``$1$'', otherwise, see
Figure~\ref{F:distance}.
%%%
%%%%%%%%%%%%%%%%%%%%%%%%%%%%%%%%%%%%%%%%%%%%%%%%%%%%%%%%%%%%%%%%%%%%%%%%%%%%%%%%
%%%
\begin{figure}
  \centering
  \includegraphics[width=\columnwidth]{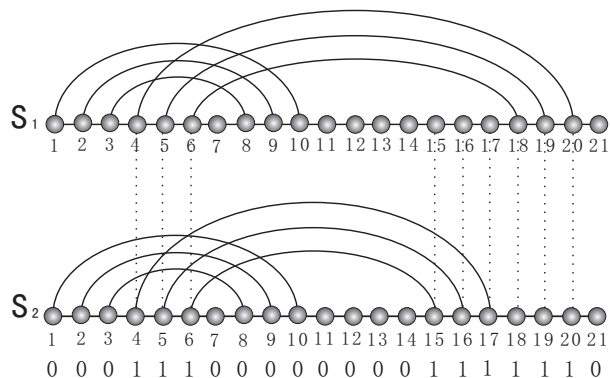}
  \caption{\small Positions paired differently in $S_1$ and $S_2$ are assigned
  a ``$1$''. There are two types of positions:
  {\bf I}. $p$ is contained in different arcs, see position $4$,
  $(4,20)\in S_1$ and $(4,17) \in S_2$.
  {\bf II}. $p$ is unpaired in one structure and
  $p$ is paired in the other, such as position $18$.
  }
  \label{F:distance}
\end{figure}
%%%
%%%%%%%%%%%%%%%%%%%%%%%%%%%%%%%%%%%%%%%%%%%%%%%%%%%%%%%%%%%%%%%%%%%%%%%%%%%%%%%
%%%

One common assumption in these inverse folding algorithms is, that the energies
of specific substructures contribute additively to the energy of the entire
structure. Let us proceed by analyzing the algorithms.

\paragraph{RNAinverse}
is the first inverse-folding algorithm that derives sequences that
realize given RNA secondary structures as mfe-configuration. In its
initialization step, a random compatible sequence $s$ for
the target $T$ is generated. Then {\tt RNAinverse} proceeds by
updating the sequence $s$ to $s',s''\dots$ step by step, minimizing
the structure distance between the mfe structure of $s'$ and the
target structure $T$. Based on the observation, that the energy of a
substructure contributes additively to the mfe of the molecule, {\tt
RNAinverse} optimizes ``small'' substructures first, eventually
extending these to the entire structure. While optimizing
substructures, RNAinverse does an adaptive walk in order to decrease
the structure distance. In fact, this walk is based entirely on
random compatible mutations.

\paragraph{RNA-SSD}
{\tt RNA-SSD} first assigns specific probabilities
to the bases located in unpaired positions and the base pairs ($\pairGC, \pairAU,
\pairUG$) of $T$, respectively.
In this assignment the probability of a unpaired position being assigned either
A or U is greater than assigning G or C. Similarly, the probability of
pairs $\pairGC$ and $\pairCG$ base pairs is greater than that of the other
base pairs.
Then, {\tt RNA-SSD} derives a hierarchical decomposition of the
target structure. It recursively splits the structure and thereby
derives a binary decomposition tree rooted in $T$ and whose leaves
correspond to $T$-substructures. Each non-leaf node of this tree
represents a substructure obtained by merging the two substructures
of its respective children. Given this tree, {\tt RNA-SSD} performs
a stochastic local search, starting at the leaves, subsequently
working its way up to the root.

\paragraph{INFO-RNA} employs a dynamic programming method for finding
a well suited initial sequence. This sequence has a lowest energy with
respect to the $T$. Since the latter does not necessarily fold into $T$,
(due to potentially existing competing configurations) {\tt INFO-RNA}
then utilizes an improved\footnote{relative to the local search
routine used in {\tt RNAinverse}} stochastic local search in order to find
a sequence in the neutral network of $T$.
In contrast to {\tt RNAinverse},
{\tt INFO-RNA} allows for increasing
the distance to the target structure. At the same time, only
positions that do not pair correctly and positions adjacent to these
are examined.

%%%
%%%%%%%%%%%%%%%%%%%%%%%%%%%%%%%%%%%%%%%%%%%%%%%%%%%%%%%%%%%%%%%%%%%%%%%%%%%%%%%%
%%%
\subsection{$\foldpk$}\label{S:loops}
%%%
%%%%%%%%%%%%%%%%%%%%%%%%%%%%%%%%%%%%%%%%%%%%%%%%%%%%%%%%%%%%%%%%%%%%%%%%%%%%%%%%
%%%

$\foldpk$ is an {\it ab initio} folding algorithm that maps RNA sequences into
$3$-noncrossing RNA structures. It is guaranteed to search all $3$-noncrossing,
$\sigma$-canonical structures and derives some (not necessarily unique),
loop-based mfe-configuration. In the following we always assume $\sigma\ge 3$.
The input of $\foldpk$ is an arbitrary RNA sequence $s$ and an integer $N$.
Its output is a list of $N$ $3$-noncrossing, $\sigma$-canonical structures,
the first of which being the mfe-structure for $s$.
This list of $N$ structures $(C_0,C_1,\dots,C_{N-1})$ is ordered by
the free energy and
the first list-element, the mfe-structure, is denoted by $\foldpk(s)$.
If no $N$ is specified, $\foldpk$ assumes $N=1$ as default.

{\tt Cross} generates a mfe-structure based on specific loop-types of
$3$-noncrossing RNA structures. For a given structure $S$, let $\alpha$ be
an arc contained in $S$ ($S$-arc) and denote the set of $S$-arcs that
cross $\alpha$ by $\mathscr{A}_{S}(\alpha)$.

For two arcs $\alpha=(i,j)$ and $\alpha'=(i',j')$, we next specify the
partial order ``$\prec$'' over the set of arcs:
\begin{equation*}
\alpha' \prec \alpha \quad \text{\rm  if and only if}\quad
i<i'<j'<j.
\end{equation*}
All notions of minimal or maximal elements are understood to be with respect to
$\prec$.
An arc $\alpha \in \mathscr{A}_S(\beta)$ is
called a minimal, $\beta$-crossing if there exists no $\alpha' \in
\mathscr{A}_S(\beta)$ such that $\alpha' \prec \alpha$. Note that $\alpha \in
\mathscr{A}_S(\beta)$ can be minimal $\beta$-crossing, while $\beta$ is not
minimal $\alpha$-crossing. $3$-noncrossing diagrams exhibit the following four
basic loop-types:

\paragraph{(1)}

A hairpin-loop is a pair
\[
  ((i,j),[i+1,j-1])
\]
where $(i,j)$ is an arc and $[i,j]$ is an interval, i.e.~a sequence of
consecutive vertices
 $(i,i+1,\dots,j-1,j)$.

\paragraph{(2)}

An interior-loop, is a sequence
\[
  ((i_1,j_1),[i_1+1,i_2-1],(i_2,j_2),[j_2+1,j_1-1]),
\]
where $(i_2,j_2)$ is nested in $(i_1,j_1)$.
That is we have
$i_1<i_2<j_2<j_1$.

\paragraph{(3)}

A multi-loop, see Figure~\ref{F:stand} \cite{Reidys:08algo}, is a sequence
\[
  ((i_1,j_1), [i_1+1,\omega_1-1], S_{\omega_{1}}^{\tau_1},
  [\tau_1+1,\omega_2-1], S_{\omega_{2}}^{\tau_2}, \dots ),
\]
where $S_{\omega_h}^{\tau_h}$ denotes a pseudoknot structure over $[\omega_h,
\tau_h]$ (i.e.~nested in $(i_1,j_1)$) and subject to the following condition:
if all $S_{\omega_h}^{\tau_h}=(\omega_h,\tau_h)$, i.e.~all substructures are
just arcs, for all $h$,  then we have $h\ge 2$).

%%%
%%%%%%%%%%%%%%%%%%%%%%%%%%%%%%%%%%%%%%%%%%%%%%%%%%%%%%%%%%%%%%%%%%%%%%%%%%%%%%%%
%%%
\begin{figure}
  \centering
  \includegraphics[width=\columnwidth]{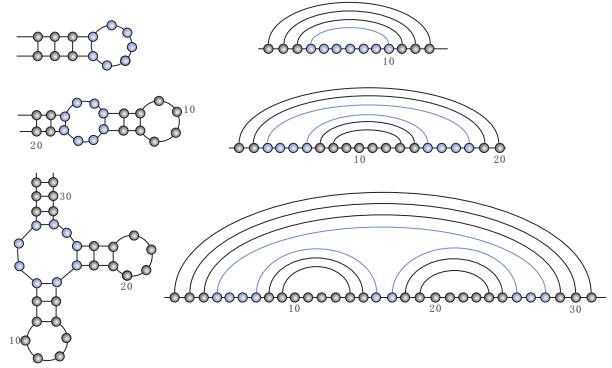}
  \caption{\small
  The standard loop-types: hairpin-loop (top), interior-loop
  (middle) and multi-loop (bottom). These represent all loop-types that occur
  in RNA secondary structures.}
  \label{F:stand}
\end{figure}
%%%
%%%%%%%%%%%%%%%%%%%%%%%%%%%%%%%%%%%%%%%%%%%%%%%%%%%%%%%%%%%%%%%%%%%%%%%%%%%%%%%%
%%%

%%%
%%%%%%%%%%%%%%%%%%%%%%%%%%%%%%%%%%%%%%%%%%%%%%%%%%%%%%%%%%%%%%%%%%%%%%%%%%%%%%%%
%%%
\begin{figure}
  \centering
  \includegraphics[width=\columnwidth]{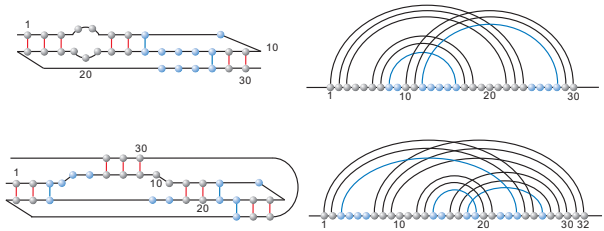}
  \caption{\small Pseudoknot loops, formed by
  all blue vertices and arcs.}
  \label{F:pseudoloop}
\end{figure}
%%%
%%%%%%%%%%%%%%%%%%%%%%%%%%%%%%%%%%%%%%%%%%%%%%%%%%%%%%%%%%%%%%%%%%%%%%%%%%%%%%%%
%%%

A pseudoknot, see Figure~\ref{F:pseudoloop} \cite{Reidys:08algo}, consists of
the following data:

\paragraph{(\sf P1)}

A set of arcs
\[
  P=\left\{(i_1,j_1),(i_2,j_2), \dots,(i_t,j_t)\right\},
\]
where $i_1=\min\{i_h\}$ and $j_t=\max\{j_h\}$, such that

\begin{enumerate}

  \item[\bf (i)] the diagram induced by the arc-set $P$ is irredu\-cible,
    i.e.~the dependency-graph of $P$ (i.e.~the graph having $P$ as vertex set
    and in which $\alpha$ and $\alpha'$ are adjacent if and only if they cross)
    is connected and

  \item[\bf (ii)] for each $(i_{h},j_{h})\in P$ there exists some arc $\beta$
    (not necessarily contained in $P$) such that $(i_{h},j_{h})$ is minimal
    $\beta$-crossing.

\end{enumerate}

\paragraph{(\sf P2)}

Any $i_1<x<j_t$, not contained in hairpin-, interior- or
multi-loops.

Having discussed the basic loop-types, we are now in position to state

%%%
%%%%%%%%%%%%%%%%%%%%%%%%%%%%%%%%%%%%%%%%%%%%%%%%%%%%%%%%%%%%%%%%%%%%%%%%%%%%%%%%
%%%
\begin{theorem}\label{th:decomp}
Any $3$-noncrossing RNA pseudoknot structure has a unique loop-decomposition
{\rm \cite{Reidys:08algo}}.
\end{theorem}

%%%
%%%%%%%%%%%%%%%%%%%%%%%%%%%%%%%%%%%%%%%%%%%%%%%%%%%%%%%%%%%%%%%%%%%%%%%%%%%%%%%%
%%%
\begin{figure}
  \centering
  \includegraphics[width=\columnwidth]{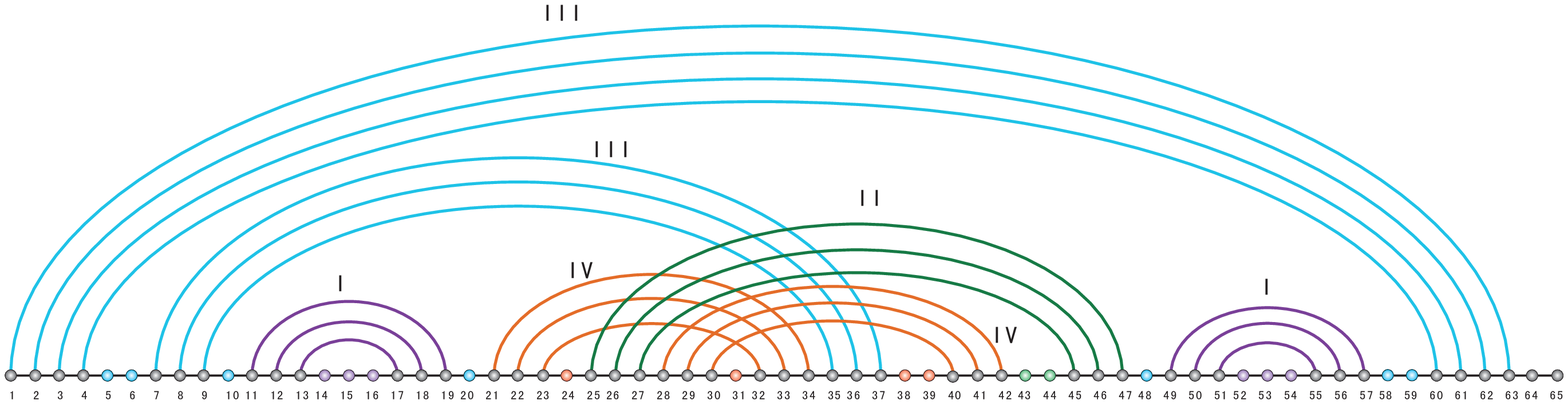}
  \caption{\small Loop decomposition: here a hairpin-loop (I),
  an interior-loop (II), a multi-loop (III) and a pseudoknot (IV).}
  \label{F:loopdecomp}
\end{figure}
%%%
%%%%%%%%%%%%%%%%%%%%%%%%%%%%%%%%%%%%%%%%%%%%%%%%%%%%%%%%%%%%%%%%%%%%%%%%%%%%%%%%
%%%

Figure~\ref{F:loopdecomp} illustrates the loop decomposition of a
$3$-noncrossing structure.

A {\bf motif} in $\foldpk$ is a $3$-noncrossing structure, having
only $\prec$-maximal stacks of size exactly $\sigma$, see Figure~\ref{F:motif}.
A {\bf skeleton}, $S$, is a ${k}$-noncrossing structure such that
\begin{itemize}
\item its core, $c(S)$ has no noncrossing arcs and
\item its $L$-graph, $L(S)$ is connected.
\end{itemize}
Here the core of a structure, $c(S)$, is obtained by collapsing its stacks into
single arcs (thereby reducing its length) and the graph $L(S)$ is obtained by
mapping arcs into vertices and connecting any two if they cross in the diagram
representation of $S$, see Figure~\ref{F:skeleton}.
As for the general strategy, $\foldpk$ constructs
$3$-noncrossing RNA structure ``from top to bottom'' via three subroutines:

%%%
%%%%%%%%%%%%%%%%%%%%%%%%%%%%%%%%%%%%%%%%%%%%%%%%%%%%%%%%%%%%%%%%%%%%%%%%%%%%%%%%
%%%
\begin{figure}
  \centering
  \includegraphics[width=\columnwidth]{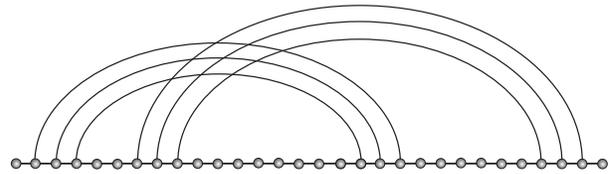}
  \caption{\small Motif: a $3$-noncrossing, $3$-canonical motif.
  }
  \label{F:motif}
\end{figure}
%%%
%%%%%%%%%%%%%%%%%%%%%%%%%%%%%%%%%%%%%%%%%%%%%%%%%%%%%%%%%%%%%%%%%%%%%%%%%%%%%%%%
%%%

%%%
%%%%%%%%%%%%%%%%%%%%%%%%%%%%%%%%%%%%%%%%%%%%%%%%%%%%%%%%%%%%%%%%%%%%%%%%%%%%%%%%
%%%
\begin{figure}
  \centering
  \includegraphics[width=\columnwidth]{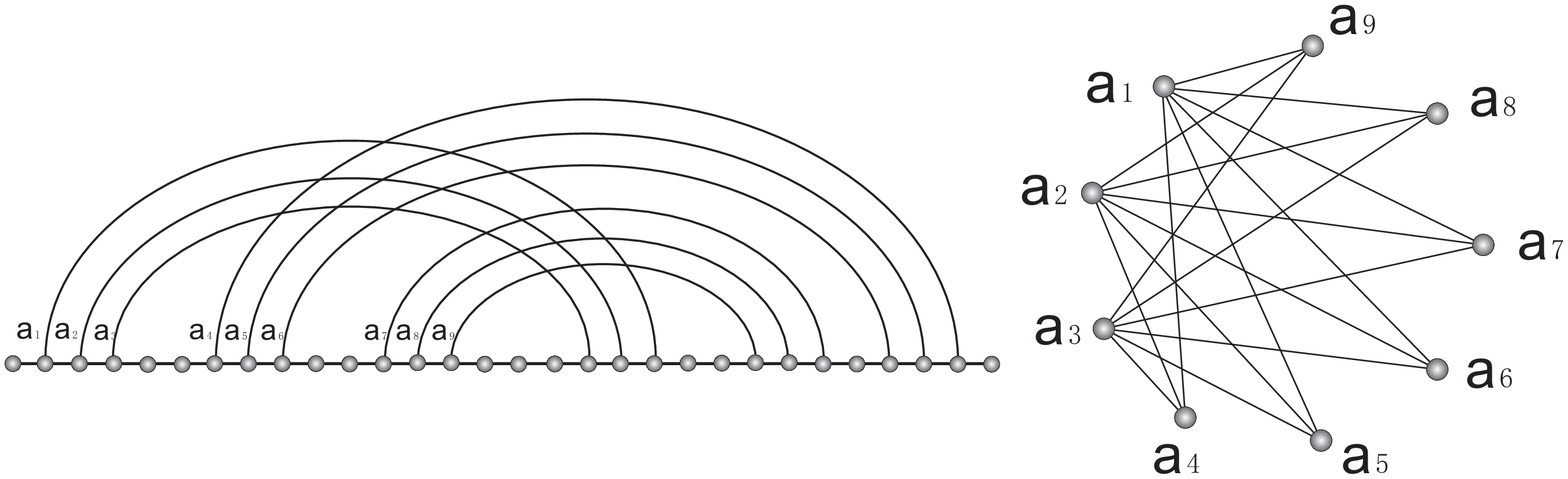}
  \caption{\small Skeleton and its
  $L$-graph: we display a skeleton (left)
  and its $L$-graph (right).
  }
  \label{F:skeleton}
\end{figure}
%%%
%%%%%%%%%%%%%%%%%%%%%%%%%%%%%%%%%%%%%%%%%%%%%%%%%%%%%%%%%%%%%%%%%%%%%%%%%%%%%%%%
%%%

%%%%%%%%%%%%%%%%%%%%%%%%%%%%%%%%%%%%%%%%%%%%%%%%%%%%%%%%%%%%%%%%%%%%%%%%%%%%%%%%
\paragraph{I ({\sc Shadow}):}

Here we generate all maximal stacks of the
structure. Note that a stack is maximal with respect to $\prec$ if it is not
nested in some other stack. This is derived by ``shadowing'' the motifs,
i.e.~their $\sigma$-stacks are extended ``from top to bottom''.
%%%%%%%%%%%%%%%%%%%%%%%%%%%%%%%%%%%%%%%%%%%%%%%%%%%%%%%%%%%%%%%%%%%%%%%%%%%%%%%%
\paragraph{II ({\sc SkeletonBranch}):}

Given a shadow, the second step of $\foldpk$ consists in
generating, the skeleta-tree. The nodes of this tree are particular
$3$-noncrossing structures, obtained by successive insertions of stacks.
Intuitively, a skeleton encapsulates all cross-serial arcs that cannot be
recursively computed. Here the tree complexity is controlled via limiting
the (total) number of pseudoknots.

%%%%%%%%%%%%%%%%%%%%%%%%%%%%%%%%%%%%%%%%%%%%%%%%%%%%%%%%%%%%%%%%%%%%%%%%%%%%%%%%
\paragraph{III ({\sc Saturation}):}

In the third subroutine each skeleton is saturated via DP-routines. After the
saturation the mfe-$3$-noncrossing structure is derived.

Figure~\ref{F:sketch} provides an overview on how the three
subroutines are combined.

%%%
%%%%%%%%%%%%%%%%%%%%%%%%%%%%%%%%%%%%%%%%%%%%%%%%%%%%%%%%%%%%%%%%%%%%%%%%%%%%%%%%
%%%
\begin{figure}
  \centering
  \includegraphics[width=\columnwidth]{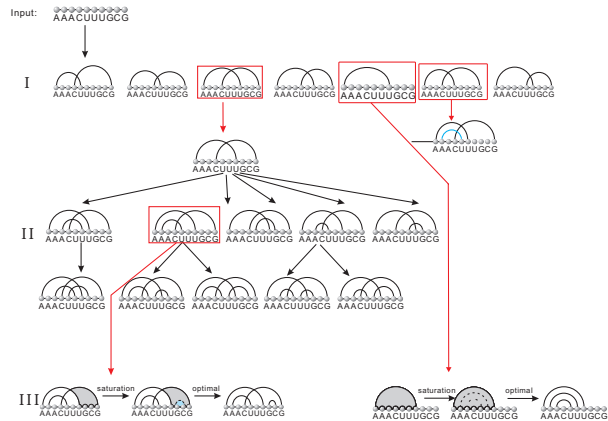}
  \caption{\small 
  An outline of $\foldpk$ (for illustration purposes we assume here $\sigma=1$):
  The routines {\sc Shadow}, {\sc SkeletonBranch} and {\sc Saturation} are
  depicted. Due to space limitations we only represent a few select motifs and
  for the same reason only one of the motifs displayed in the first row
  is extended by one arc (drawn in blue).
  Furthermore note that only motifs with crossings give rise to nontrivial
  skeleton-trees, all other motifs are considered directly as input for
  {\sc Saturation}.}
  \label{F:sketch}
\end{figure}
%%%
%%%%%%%%%%%%%%%%%%%%%%%%%%%%%%%%%%%%%%%%%%%%%%%%%%%%%%%%%%%%%%%%%%%%%%%%%%%%%%%%
%%%

%%%
%%%%%%%%%%%%%%%%%%%%%%%%%%%%%%%%%%%%%%%%%%%%%%%%%%%%%%%%%%%%%%%%%%%%%%%%%%%%%%%%
%%%
\section{The algorithm}\label{S:algo}
%%%
%%%%%%%%%%%%%%%%%%%%%%%%%%%%%%%%%%%%%%%%%%%%%%%%%%%%%%%%%%%%%%%%%%%%%%%%%%%%%%%%
%%%
The inverse folding algorithm {\tt Inv} is based on the {\it ab
initio} folding algorithm $\foldpk$. The input of ${\tt Inv}$ is the target
structure, $T$. The latter is expressed as a character string of
``\verb|:()[]{}|'', where ``{\tt :}'' denotes unpaired base and
``\verb|()|'', ``\verb|[]|'', ``\verb|{}|'' denote paired bases.

In Algorithm~\ref{A:invfold}, we present the pseudocodes of
algorithm ${\tt Inv}$. After validation of the target structure
(lines $2$~to $5$ in Algorithm~\ref{A:invfold}), similar to {\tt
INFO-RNA}, {\tt Inv} constructs an initial sequence  and
then proceeds by a stochastic local search based on the loop
decomposition of the target. This sequence is derived via the
routine $\adjustseq$. We then decompose the target structure into
loops and endow these with a linear order. According to this order
we use the routine $\aw$ in order to find for each loop a ``proper''
local solution.

%%%%%%%%%%%%%%%%%%%%%%%%%%%%%%%%%%%%%%%%%%%%%%%%%%%%%%%%%%%%%%%%%%%%%%%%%%%%%%%%
%%%
\begin{algorithm}\small
  \caption{{\tt Inv}}
  \label{A:invfold}
  \begin{algorithmic}[1]
  \item[\textbf{Input:}] $k$-noncrossing target structure $T$
  \item[\textbf{Output:}] an RNA sequence $seq$
    \REQUIRE $k \le 3$ and $T$ is composed with ``\verb|:()[]{}|''
    \ENSURE $\foldpk(seq) = T$

    \STATECOMMENT Step 1: Validate structure
    \IF{${\bf false} = \checkstru(T)$}
      \MYPRINT incorrect structure
      \MYRETURN \texttt{NIL}
    \ENDIF

    \STATE
    \STATECOMMENT{Step 2: Generate the start sequence}
    \STATE $start \leftarrow \makestart(T)$

    \STATE
    \STATECOMMENT{Step 3: Adjust the start sequence}
    \STATE $seq_\textrm{middle} \leftarrow \adjustseq(start, T)$

    \STATE
    \STATECOMMENT{Step 4: Decompose $T$ and derive the ordered intervals.}
    \STATE Interval array $I$
    \STATE $m \leftarrow \vert I\vert$
    \COMMENT{$I$ satisfies $I_m=T$}

    \STATE
    \STATECOMMENT{Step 5: Stochastic Local Search}
    \STATE $seq\leftarrow seq_\textrm{middle}$
    \FORALL{intervals in the array $I_w$}
      \STATE $l \leftarrow \text{start-point}(I_w)$
      \STATE $r \leftarrow \text{end-point}(I_w)$
      \STATE $s' \leftarrow seq|_{[l,r]}$
      \COMMENT{get sub-sequence}
      \STATE $seq|_{[l,r]}\leftarrow \aw(s',I_w)$
    \ENDFOR

    \STATE
    \STATECOMMENT{Step 6: output}
    \IF{$seq_\textrm{min} = \foldpk(seq)$}
      \MYRETURN $seq$
    \ELSE
      \MYPRINT Failed!
      \MYRETURN {\tt NIL}
    \ENDIF
  \end{algorithmic}
\end{algorithm}
%%%
%%%%%%%%%%%%%%%%%%%%%%%%%%%%%%%%%%%%%%%%%%%%%%%%%%%%%%%%%%%%%%%%%%%%%%%%%%%%%%%%
%%%

\subsection{{\sc Adjust-Seq}}\label{S:init}

In this section we describe Steps $2$~and $3$ of the pseudocodes presented in
Algorithm~\ref{A:invfold}. The routine $\makestart$, see line $8$, generates a
random sequence, $start$, which is compatible to the target, with uniform
probability.

We then initialize the variable $seq_\textrm{min}$ via the sequence $start$ and
set the variable $d=+\infty$, where $d$ denotes the structure
distance between $\foldpk(seq_\textrm{min})$
and $T$.

Given the sequence $start$, we construct a set of potential ``competitors'',
$C$, i.e.~a set of structures suited as folding targets for $start$. In
Algorithm~\ref{A:adjust} we show how to adjust the start sequence using the
routine $\adjustseq$. Lines $4$~to $38$ of Algorithm~\ref{A:adjust}, contain a
{\bf For}-loop, executed at most $\sqrt{n}/2$ times. Here the loop-length
$\sqrt{n}/2$ is heuristically determined.

Setting the $\foldpk$-parameter\footnote{For all computer experiments
we set $N=50$.}, $N$,
 the subroutine executed in the
loop-body consists of the following three steps.

%%%%%%%%%%%%%%%%%%%%%%%%%%%%%%%%%%%%%%%%%%%%%%%%%%%%%%%%%%%%%%%%%%%%%%%%%%%%%%%%
\paragraph{Step I. Generating $C^0(\lambda^i)$ via $\foldpk$.}

Suppose we are in the $i$th step of the {\bf For}-loop and are given
the sequence $\lambda^{i-1}$ where $\lambda^0=start$. We consider
$\foldpk(\lambda^{i-1},N)$, i.e.~the list of suboptimal structures
with respect to $\lambda^{i-1}$,
$$C^0(\lambda^{i-1})=\foldpk(\lambda^{i-1},N)=(C_h^{0}(\lambda^{i-1}))_{h=0}^{N-1} $$
If $C_0^0(\lambda^{i-1})=T$, then {\tt
Inv} returns $\lambda^{i-1}$. Else, in case of
$d=(\foldpk(C_0^0(\lambda^{i-1})),T)<d_{min}$, we set
\begin{eqnarray*}
seq_\textrm{min} &= & \lambda^{i-1} \\
d_\textrm{min} &= & d(\foldpk(C_0^0(\lambda^{i-1})),T).
\end{eqnarray*}
Otherwise we do not update
$seq_\textrm{min}$ and go directly to Step II.

%%%%%%%%%%%%%%%%%%%%%%%%%%%%%%%%%%%%%%%%%%%%%%%%%%%%%%%%%%%%%%%%%%%%%%%%%%%%%%%%
\paragraph{Step II. The competitors.}

We introduce a specific procedure that
``perturbs'' arcs of a given RNA pseudoknot structure, $S$. Let $a$ be an arc
of $S$ and let $l(a)$, $r(a)$ denote the start- and end-point of $a$. A
perturbation of $a$ is a procedure which generates a new arc $a'$, such that
\[
  \vert l(a)-l(a')\vert\le1 \quad {\rm and} \quad \vert r(a)-r(a')\vert\le1\,.
\]
Clearly, there are nine perturbations of any given arc $a$ (including $a$
itself), see Figure~\ref{F:adjacentarc}.

We proceed by keeping $a$, replacing the arc $a$ by a nontrivial perturbation
or remove $a$, arriving at a set of ten structures $\nu(S,a)$.

%%%
%%%%%%%%%%%%%%%%%%%%%%%%%%%%%%%%%%%%%%%%%%%%%%%%%%%%%%%%%%%%%%%%%%%%%%%%%%%%%%%%
%%%
\begin{figure}
  \centering
  \includegraphics[width=\columnwidth]{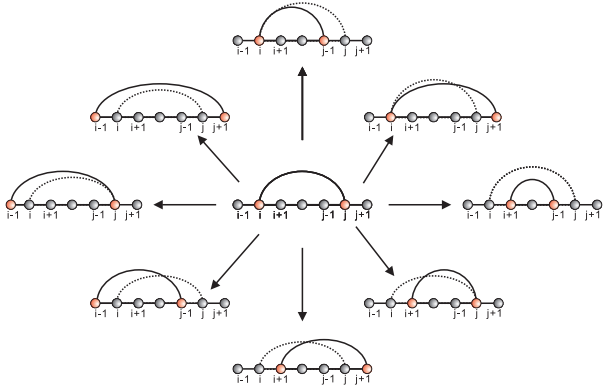}
  \caption{\small Nine perturbations of an arc $(i,j)$. Original arcs are drawn
  dotted, and the arcs incident to red bases are the perturbations.}
  \label{F:adjacentarc}
\end{figure}
%%%
%%%%%%%%%%%%%%%%%%%%%%%%%%%%%%%%%%%%%%%%%%%%%%%%%%%%%%%%%%%%%%%%%%%%%%%%%%%%%%%%
%%%

Now we use this method in order to generate the set $C^1(\lambda^{i-1})$ by
perturbing each arc of each structure $C^0_h(\lambda^{i-1})\in
C^0(\lambda^{i-1})$.
If $C^0_h(\lambda^{i-1})$ has $A_h$ arcs, $\{a_h^1, \dots, a_h^{A_h}\}$, then
\[
  C^1(\lambda^{i-1})=\bigcup_{h=0}^{N-1}\bigcup_{j=1}^{\;A_h}
\nu(C^0_h(\lambda^{i-1}),a_h^j)\,.
\]
This construction may result in duplicate, inconsistent or incompatible
structures. Here, a structure is inconsistent if there exists at least one
position paired with more than one base, and incompatible if there exists at
least one arc not compatible with $\lambda^{i-1}$, see Figures
\ref{F:inconsistent.arc}~and \ref{F:incompatible.arc}. Here compatibility is
understood with respect to the Watson-Crick and ${\bf G}$-${\bf U}$ base pairing
rules. Deleting inconsistent and incompatible structures, as well as those identical
to the target, we arrive at the set of competitors,
\begin{equation*}
C(\lambda^{i-1})\subset C^1(\lambda^{i-1}).
\end{equation*}

%%%
%%%%%%%%%%%%%%%%%%%%%%%%%%%%%%%%%%%%%%%%%%%%%%%%%%%%%%%%%%%%%%%%%%%%%%%%%%%%%%%
%%%

\begin{figure}
  \centering
  \includegraphics[width=0.7\columnwidth]{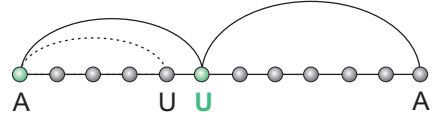}
  \caption{\small Inconsistent structures: the dotted arc is perturbed by
  shifting its end-point. This perturbation leads to a nucleotide establishing
  two base pairs, which is impossible.}
  \label{F:inconsistent.arc}
\end{figure}
%%%
%%%%%%%%%%%%%%%%%%%%%%%%%%%%%%%%%%%%%%%%%%%%%%%%%%%%%%%%%%%%%%%%%%%%%%%%%%%%%%%
%%%
\begin{figure}
  \centering
  \includegraphics[width=0.7\columnwidth]{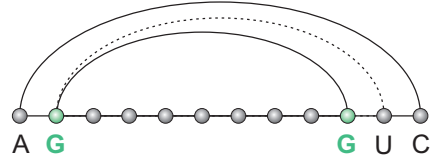}
  \caption{\small
  Incompatible structures: we display a perturbation of the dotted arc
  leading to a structure that is incompatible to the given sequence.}
  \label{F:incompatible.arc}
\end{figure}
%%%
%%%%%%%%%%%%%%%%%%%%%%%%%%%%%%%%%%%%%%%%%%%%%%%%%%%%%%%%%%%%%%%%%%%%%%%%%%%%%%%
%%%

%%%
%%%%%%%%%%%%%%%%%%%%%%%%%%%%%%%%%%%%%%%%%%%%%%%%%%%%%%%%%%%%%%%%%%%%%%%%%%%%%%%
%%%
\begin{figure}
  \centering
  \includegraphics[width=\columnwidth]{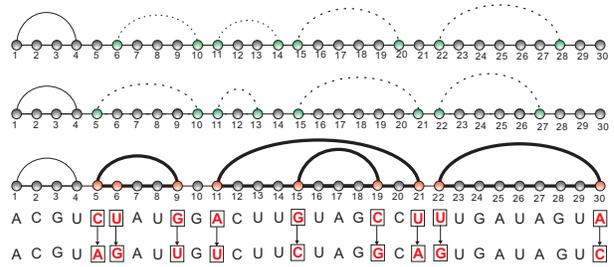}
  \caption{\small
  Mutation: Suppose the top and middle structures represent the set of competitors
  and the bottom structure is target. We display $\lambda^{i-1}$
  (top sequence) and its mutation, $\lambda^{i}$ (bottom sequence).
  Two nucleotides of base pairs not contained in $T$ are colored green, nucleotides
  subject to mutations are colored red.}
  \label{F:adjustseq}
\end{figure}
%%%
%%%%%%%%%%%%%%%%%%%%%%%%%%%%%%%%%%%%%%%%%%%%%%%%%%%%%%%%%%%%%%%%%%%%%%%%%%%%%%%
%%%

\paragraph{Step III. Mutation}

Here we adjust $\lambda^{i-1}$ with respect to $T$ as well as the set of
competitors, $C(\lambda^{i-1})$ derived in the previous step.
Suppose $\lambda^{i-1}=s_1^{i-1} s_2^{i-1}
\dots s_n^{i-1}$. Let
$p(S,w)$ be the position paired to the position $w$ in the RNA structure $S\in
C(\lambda^{i-1})$, or 0 if position $w$ is unpaired. For instance, in
Figure~\ref{F:adjustseq}, we have $p(T,1)=4$, $p(T,2)=0$ and $p(T,4)=1$. For
each position $w$ of the target $T$, if there exists a structure
$C_h(\lambda^{i-1})\in C(\lambda^{i-1})$ such that $p(C_h(\lambda^{i-1}),w)
\not=p(T,w)$
(see positions $5$, $6$, $9$, and $11$ in Figure~\ref{F:adjustseq}) we modify
$\lambda^{i-1}$ as follows:

\begin{enumerate}
  \item {\bf unpaired position:} If $p(T,w)=0$, we update
	  $s^{i-1}_w$ randomly into the nucleotide $s^{i}_w\neq s^{i-1}_w$,
	  such that for each $C_h(\lambda^{i-1})\in C(\lambda^{i-1})$,
	  either $p(C_h(\lambda^{i-1}),w)=0$ or $s^{i}_w$ is
	  not compatible with $s^{i-1}_v$
      where $v=p(C_h(\lambda^{i-1}),w)>0$, See
      position $6$ in Figure~\ref{F:adjustseq}.

  \item {\bf start-point:} If $p(T,w)>w$, set $v=p(T,w)$. We randomly choose a
	  compatible base pair $(s^{i}_w,s^{i}_v)$ different
	  from $(s^{i-1}_w,s^{i-1}_v)$, such
      that for each $C_h(\lambda^{i-1})\in C(\lambda^{i-1})$, either
      $p(C_h(\lambda^{i-1}),w)=0$ or $s^{i}_w$ is not compatible with
      $s^{i-1}_u$, where
      $u=p(C_h(\lambda^{i-1}),w)>0$ is the end-point paired with $s_w^{i-1}$ in
      $C_h(\lambda^{i-1})$ (Figure~\ref{F:adjustseq}: $(5,9)$. The pair
      \pairGC{} retains the compatibility to $(5,9)$, but is incompatible
      to $(5,10)$). By Figure~\ref{F:adjustproof} we show feasibility of this step.

  \item {\bf end-point:} If $0<p(T,w)<w$, then by construction the nucleotide has
	already been considered in the previous step.
\end{enumerate}

%%%
%%%%%%%%%%%%%%%%%%%%%%%%%%%%%%%%%%%%%%%%%%%%%%%%%%%%%%%%%%%%%%%%%%%%%%%%%%%%%%
%%%
\begin{figure}
  \centering
  \includegraphics[width=\columnwidth]{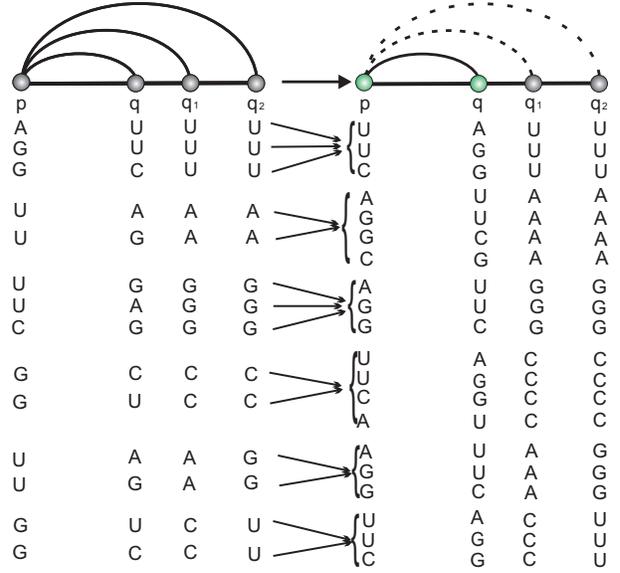}
  \caption{\small Mutations are always possible:
  suppose $p$ is paired with $q$ in $T$
  and $p$ is paired with $q_1$ in one competitor and $q_2$ in another one.
  For a fixed nucleotide at $p$ there are at most two scenarios, since a
  base can pair with at most two different bases.
  For instance, for $\baseG$ we have the pairs ${\pairGC,\pairGU}$.
  We display all nucleotide configurations (LHS) and their
  corresponding solutions (RHS).}
  \label{F:adjustproof}
\end{figure}
%%%
%%%%%%%%%%%%%%%%%%%%%%%%%%%%%%%%%%%%%%%%%%%%%%%%%%%%%%%%%%%%%%%%%%%%%%%%%%%%%%
%%%

Therefore, updating all the nucleotides of $\lambda^{i-1}$, we arrive at the
new sequence $\lambda^i=s^{i}_{1}s^{i}_2\dots s^i_n$.

Note that the above mutation steps
heuristically decrease the structure distance. However, the resulting sequence is
not necessarily incompatible to all competitors. For instance, consider a competitor
$C_h$ whose arcs are all contained $T$.
Since $\lambda^i$ is compatible with $T$,
$\lambda^i$ is compatible with $C_h$. Since competitors
are obtained from suboptimal folds such a scenario may arise.

In practice, this situation represents not a problem, since these type of
competitors are likely to be ruled out by virtue of the fact that they
have a mfe larger than that of the target structure.

Accordingly we have the following situation, competitors are eliminated due
to two, equally important criteria: incompatibility as well as minimum free
energy considerations.

If the distance of $\foldpk(\lambda^{i})$ to $T$ is less
than or equal to $d_\textrm{min}+5$, we return to Step I (with
$\lambda^{i}$).  Otherwise, we repeat Step III (for at most 5 times) thereby
generating $ \lambda^{i}_1,\dots,\lambda_5^{i}$ and set
$\lambda^{i}=\lambda^{i}_w$ where $d(\foldpk(\lambda_w^{i}),T)$
is minimal.

The procedure $\adjustseq$ employs the negative paradigm \cite{Dirks:2004}
in order to exclude energetically close conformations. It returns the sequence
$seq_\textrm{middle}$ which is tailored to realize the target structure as
mfe-fold.

%%%
%%%%%%%%%%%%%%%%%%%%%%%%%%%%%%%%%%%%%%%%%%%%%%%%%%%%%%%%%%%%%%%%%%%%%%%%%%%%%%%
%%%
\begin{algorithm}\small
  \caption{$\adjustseq$}\label{A:adjust}
  \begin{algorithmic}[1]
  \item[\textbf{Input:}] the original start sequence $start$
  \item[\textbf{Input:}] the target structure $T$
  \item[\textbf{Output:}] a initialized sequence $seq_\textrm{middle}$
    \STATE $n \leftarrow \textrm{length of } T$
    \STATE $d_\textrm{min} \leftarrow +\infty$,
      \quad $seq_\textrm{min} \leftarrow start$
    \FOR{$i=1$ to $\frac{1}{2}\sqrt{n}$}
      \STATECOMMENT{Step I: generate the set $C^0(\lambda^{i-1})$ via
      $\foldpk$}
      \STATE $C^0(\lambda^{i-1}) \leftarrow \foldpk(\lambda^{i-1},N)$

      \STATE $d\leftarrow d(C_0^0(\lambda^{i-1}), T)$
      \IF{$d=0$}
        \MYRETURN $\lambda^{i-1}$
      \ELSIF{$d<d_\textrm{min}$}
        \STATE $d_\textrm{min}\leftarrow d$,
          \quad $seq_\textrm{min}\leftarrow \lambda^{i-1}$
      \ENDIF

      \STATE
      \STATECOMMENT{Step II: generate the competitor set $C(\lambda^{i-1})$}
      \STATE $C^1(\lambda^{i-1}) \leftarrow \phi$
      \FORALL{$C_h^1(\lambda^{i-1}) \in C^1(\lambda^{i-1})$}
        \FORALL{arc $a_h^j$ of $C_h^1(\lambda^{i-1})$}
          \STATE $C^1(\lambda^{i-1})\leftarrow C^1(\lambda^{i-1})
          \cup\nu(C_0^1(\lambda^i),a_h^j)$
        \ENDFOR
      \ENDFOR
      \STATE $C(\lambda^{i-1})=$
      \STATE $\{C^1_h(\lambda^{i-1})\in C^1(\lambda^{i-1})
       \colon C^1_h(\lambda^{i-1}) \textrm{is valid}\}$

      \STATE
      \STATECOMMENT{Step III: mutation}
      \STATE $seq \leftarrow \lambda^{i-1}$
      \FOR{$w=1$ to $n$}
        \IF{$\exists C_h(\lambda^{i-1})\in C(\lambda^{i-1}) $ s.t.
       $p(C_h,w)\not=p(T,w)$}
          \STATE $seq[w] \leftarrow$ random nucleotide or pair, s.t.
            $\forall C_h(\lambda^{i-1})\in C(\lambda^{i-1})$, $seq\in C[T]$ and
       $seq\notin C[C_h(\lambda^{i-1})]$.
        %% NOTE!!!
        %% here 'C[S]' means the set of all compatible sequences of
       %% structure S
        \ENDIF
      \ENDFOR

      \STATE $T_{seq}\leftarrow \foldpk(seq)$
      \IF{$d(T_{seq}, T)<d_\textrm{min}+5$}
        \STATE $seq_\textrm{middle}\leftarrow seq$
      \ELSIF{Step III run less than 5 times}
        \GOTO Step III
      \ENDIF

    \ENDFOR\COMMENT{loop to line 3}

    \STATE
    \MYRETURN $seq_\textrm{middle}$
  \end{algorithmic}
\end{algorithm}
%%%
%%%%%%%%%%%%%%%%%%%%%%%%%%%%%%%%%%%%%%%%%%%%%%%%%%%%%%%%%%%%%%%%%%%%%%%%%%%%%%%
%%%

\subsection{{\sc Decompose} and {\sc Local-Search}}\label{S:sls}

In this section we introduce two the routines, {\sc Decompose} and
{\sc Local-Search}. The routine $\decompose$ partitions $T$ into
linearly ordered energy independent components, see
Figure~\ref{F:loopdecomp} and Section~\ref{S:loops}. {\sc
Local-Search} constructs iteratively an optimal sequence for $T$ via
local solutions, that are optimal to certain substructures of $T$.

$\decompose$: Suppose $T$ is decomposed as follows,
\[
  B=\{T_1, \dots, T_{m'}\}\,.
\]
where the $T_w$ are the loops together with all arcs in the
associated stems of the target.

 We define a linear order over $B$ as
follows: $T_w<T_h$ if either
\begin{enumerate}
  \item $T_w$ is nested in $T_h$, or
  \item the start-point of $T_w$ precedes that of $T_h$.
\end{enumerate}

In Figure~\ref{F:looporder} we display the linear order of the loops of the
structure shown in Figure~\ref{F:loopdecomp}.

%%%
%%%%%%%%%%%%%%%%%%%%%%%%%%%%%%%%%%%%%%%%%%%%%%%%%%%%%%%%%%%%%%%%%%%%%%%%%%%%%%%
%%%
\begin{figure}
  \centering
  \includegraphics[width=\columnwidth]{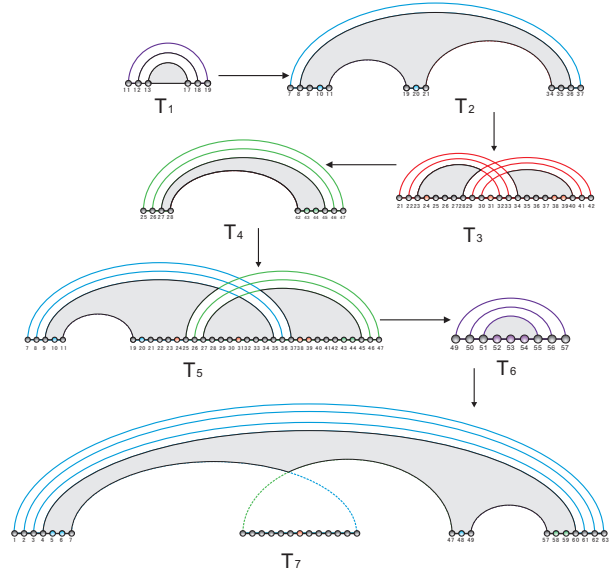}
  \caption{\small Linear ordering of loops:
  $a_1=[11,19]$, $b_1=[10,20]$, $a_2=[7,37]$, $b_2=[5,39]$,
  $a_3=[21,42]$, $b_3=[20,44]$, $a_4=[25,47]$,
  $b_4=[24,48]$, $a_5=[7,47]$, $b_5=[5,48]$, $a_6=[49,57]$,
  $b_6=[48,59]$, $a_7=[1,63]$, $b_7=[1,65]$.}
  \label{F:looporder}
\end{figure}
%%%
%%%%%%%%%%%%%%%%%%%%%%%%%%%%%%%%%%%%%%%%%%%%%%%%%%%%%%%%%%%%%%%%%%%%%%%%%%%%%%%
%%%

Next we define the interval
\begin{equation*}
a_w=[l(T_w), r(T_w)]\quad 1\leq w \leq m',
\end{equation*}
projecting the loop
$T_w$ onto the interval $[l(T_w), r(T_w)]$ and $b_w= [l', r']
\supset a_w$, being the maximal interval consisting of $a_w$ and its
adjacent unpaired consecutive nucleotides, see Figure~\ref{F:loopdecomp}.
Given two consecutive loops
$T_w <T_{w+1}$, we have two scenarios:
\begin{itemize}
\item either $b_w$ and $b_{w+1}$ are adjacent, see $b_5$ and $b_6$ in
Figure~\ref{F:looporder},
\item or $b_w \subseteq b_{w+1}$, see $b_1$ and $b_2$ in
Figure~\ref{F:looporder}.
\end{itemize}
Let $c_w=\cup_{h=1}^{w}b_h$, then
we have the sequence of intervals $a_1,b_1,c_1,\dots,a_{m'},b_{m'},c_{m'}$.
If there are no unpaired nucleotides adjacent to $a_w$, then $a_w=b_w$
and we simply delete all such $b_w$.
Thereby we derive the sequence of intervals
$I_1,I_2,\dots,I_{m}$. In Figure~\ref{F:exam.decomp.2} we illustrate
how to obtain this interval sequence: here the target
decomposes into the loops $T_1$, $T_2$ and we have $I_1=[3,5]$, $I_2=
[3,6]$, $I_3 = [2,9]$, and $I_4=[1,10]$.

%%%
%%%%%%%%%%%%%%%%%%%%%%%%%%%%%%%%%%%%%%%%%%%%%%%%%%%%%%%%%%%%%%%%%%%%%%%%%%%%%%%
%%%
\begin{figure}
  \centering
  \includegraphics[width=0.7\columnwidth]{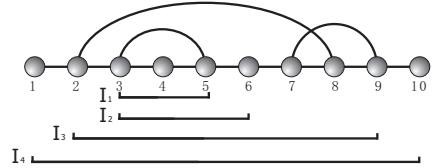}
  \caption{\small Loops and their induced sequence of intervals.}
  \label{F:exam.decomp.2}
\end{figure}
%%%
%%%%%%%%%%%%%%%%%%%%%%%%%%%%%%%%%%%%%%%%%%%%%%%%%%%%%%%%%%%%%%%%%%%%%%%%%%%%%%%
%%%

$\aw$: Given the sequence of intervals $I_1,I_2,\cdots,I_m$.
We proceed by performing a local stochastic search
on the subsequences $seq|_{I_{1}}, seq|_{I_{2}},\dots, seq|_{I_{m}}$
(initialized via $seq=seq_\textrm{middle}$ and where $s|_{[x,y]}=s_xs_{x+1}
\dots s_y$).
When we perform the local search on $seq|_{I_w}$, only positions that
contribute to the distance to the target, see Figure~\ref{F:distance},
or positions adjacent to the latter, will
be altered. We use the arrays $U_1$, $U_2$ to store
the unpaired and paired positions of $T$.
In this process, we allow for mutations that increase the structure
distance by five with probability $0.1$. 
The latter parameter is heuristically determined.
We iterate this routine until the distance is either zero or some
halting criterion is met.

%%%
%%%%%%%%%%%%%%%%%%%%%%%%%%%%%%%%%%%%%%%%%%%%%%%%%%%%%%%%%%%%%%%%%%%%%%%%%%%%%%%
%%%
\begin{algorithm}\small
  \caption{$\aw$}\label{A:sls}
  \begin{algorithmic}[1]
  \item[\textbf{Input:}] $seq_\textrm{middle}$
  \item[\textbf{Input:}] the target $T$
  \item[\textbf{Output:}] $seq$
  \ENSURE $\foldpk(seq) = T$
  \STATE $seq\leftarrow seq_\textrm{middle}$
   \IF{$\foldpk(seq)=T$}
      \MYRETURN $seq$
   \ENDIF
  \STATE decompose $T$ and derive the ordered intervals.
  \STATE $I\leftarrow[I_1,I_2,\dots,I_m]$
     \FORALL{$I_w$ in $I$}
       \STATECOMMENT{Phase I: Identify positions.}

       \STATE $d_{min}=d(\foldpk(seq|_{I_w},T|_{I_w})$
       \COMMENT{initialize $d_\textrm{min}$}

       \STATE
       \STATE derive $U_1$ via $\foldpk(seq|_{I_w})$,$T|_{I_w}$
       \STATE derive $U_2$ via $\foldpk(seq|_{I_w})$,$T|_{I_w}$
       \STATE
       \STATECOMMENT{Phase II: Test and Update.}
       \FORALL {$p$ in $U_1$}
         \STATE random $T$ compatible mutate $seq_p$
       \ENDFOR
       \FORALL {$[p,q]$ in $U_2$}
         \STATE random $T$ compatible mutate $seq_p$
       \ENDFOR

       \STATE
       \STATE $E\leftarrow\phi$
       \FORALL{$p\in U_1, U_2$}
         \STATE
	 \STATE
         \STATE $d\leftarrow d(T,\foldpk(seq_p))$
         \IF{$d<d_{min}$}
           \STATE $d_\textrm{min}\leftarrow d, \quad seq\leftarrow seq_p$
           \GOTO Phase I
         \ELSIF{$d_{min}<d<d_{min}+5$}
           \GOTO Phase I with the probability $0.1$
         \ENDIF
        \IF{$d=d_{min}$}
           \STATE $E\leftarrow E\cup\{seq\}$
        \ENDIF
      \ENDFOR

      \STATE $seq\leftarrow e_0\in E$, where $e_0$ has the lowest mfe in $E$
      \IF{Phase I run less than $10\,n$ times}
        \GOTO Phase I
      \ENDIF
      \ENDFOR
     \MYRETURN $seq$
  \end{algorithmic}
\end{algorithm}

%%%
%%%%%%%%%%%%%%%%%%%%%%%%%%%%%%%%%%%%%%%%%%%%%%%%%%%%%%%%%%%%%%%%%%%%%%%%%%%%%%%
%%%

\section{Discussion}\label{S:results}
%%%
%%%%%%%%%%%%%%%%%%%%%%%%%%%%%%%%%%%%%%%%%%%%%%%%%%%%%%%%%%%%%%%%%%%%%%%%%%%%%%%
%%%
The main result of this paper is the presentation of the algorithm
{\tt Inv}, freely available at
\begin{center}
  \url{http://www.combinatorics.cn/cbpc/inv.html}
\end{center}
Its input is a $3$-noncrossing RNA structure $T$, given in terms of
its base pairs $(i_1, i_2)$ (where $i_1<i_2$). The output of {\tt
Inv} is an RNA sequences $s=(s_1 s_2 \dots s_n)$, where $s_h \in
\{\baseA, \baseC, \baseG, \baseG\}$ with the property $\foldpk(s) =
T$, see Figure~\ref{F:UTR_sequence}.
%%%
%%%%%%%%%%%%%%%%%%%%%%%%%%%%%%%%%%%%%%%%%%%%%%%%%%%%%%%%%%%%%%%%%%%%%%%%%%%%%%%
%%%
\begin{figure}
  \centering
  \includegraphics[width=\columnwidth]{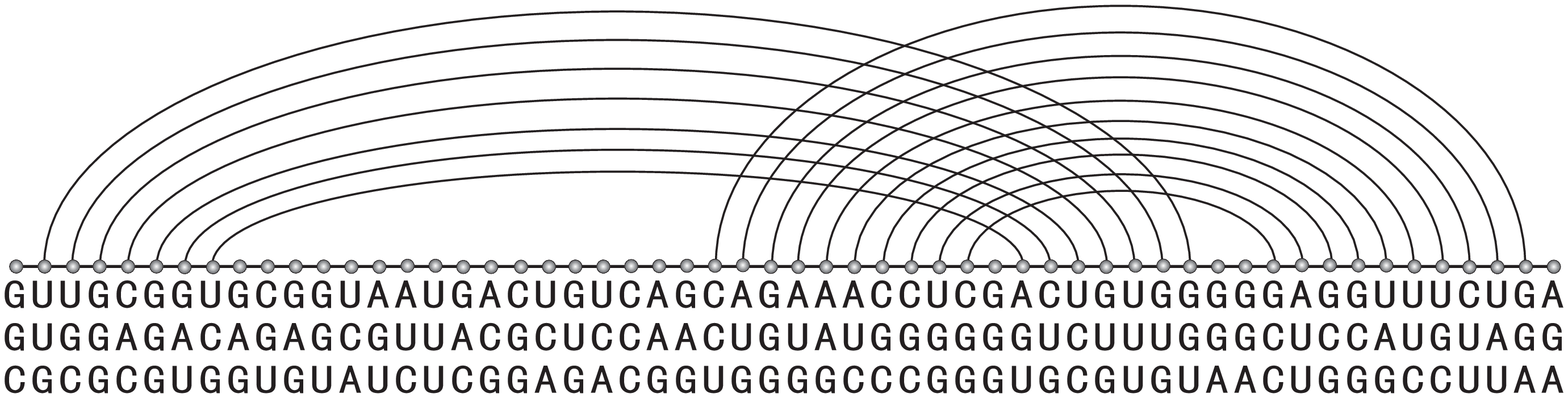}
  \caption{\small UTR pseudoknot of bovine coronavirus \cite{W:UTR}:
  its diagram representation and three sequences of its neutral network as
  constructed by {\tt Inv}.}
  \label{F:UTR_sequence}
\end{figure}

The core of {\tt Inv} is a stochastic local search
routine which is based on the fact that each $3$-noncrossing RNA
structure has a unique loop-decomposition, see
Theorem~\ref{th:decomp} in Section~\ref{S:loops}. {\tt Inv}
generates ``optimal'' subsequences and eventually arrives at a
global solution for $T$ itself.
{\tt Inv} generalizes the existing inverse folding algorithm by
considering arbitrary $3$-noncrossing canonical pseudoknot
structures. Conceptually, {\tt Inv} differs from {\tt INFO-RNA} in
how the start sequence is being generated and the particulars of the
local search itself.

As discussed in the introduction it has to be given an argument as to
{\it why} the inverse folding of pseudoknot RNA structures works.
While folding maps into RNA secondary structures are
well understood, the generalization to $3$-noncrossing RNA structures
is nontrivial. However the combinatorics of RNA pseudoknot structures
\cite{Reidys:07pseu,Reidys:08lego,Reidys:08central} implies the existence
of large neutral networks, i.e.~networks composed by sequences that all
fold into a specific pseudoknot structure. Therefore, the fact that
it is indeed possible to generate via {\tt Inv} sequences contained in the
neutral networks of targets against competing pseudoknot
configurations, see Figure~\ref{F:UTR_sequence} and
Figure~\ref{F:CrPV_IRES-PKI} confirms the predictions of \cite{Reidys:08ma}.
%%%
%%%%%%%%%%%%%%%%%%%%%%%%%%%%%%%%%%%%%%%%%%%%%%%%%%%%%%%%%%%%%%%%%%%%%%%%%%%%%%%
%%%
\begin{figure}
  \centering
  \includegraphics[width=\columnwidth]{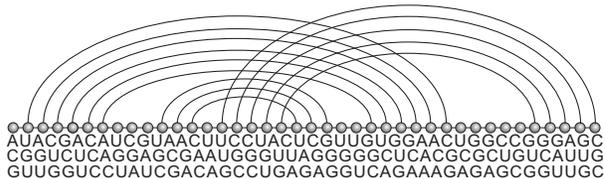}
  \caption{\small The Pseudoknot PKI of the internal ribosomal entry
  site (IRES) region \cite{W:Crpv}:
  its diagram representation and three sequences of its neutral network as
  constructed by {\tt Inv}.}
  \label{F:CrPV_IRES-PKI}
\end{figure}
%%%

An interesting class are the $3$-noncrossing nonplanar pseudoknot structures.
A nonplanar pseudoknot structure
is a $3$-noncrossing structure which is not a bi-secondary structure in the
sense of Stadler \cite{Stadler:99}. That is, it cannot be represented by
noncrossing arcs using the upper and lower half planes.
Since DP-folding paradigms of pseudoknots folding are based on gap-matrices
\cite{RE:98}, the minimal class of ``missed'' structures\footnote{given
the implemented truncations} are exactly these, nonplanar, $3$-noncrossing
structures. In Figure~\ref{F:pk} we showcase a nonplanar RNA
pseudoknot structure and $3$ sequences of its neutral network, generated by
{\tt Inv}.

As for the complexity of {\tt Inv}, the determining factor is the
subroutine $\aw$. Suppose that the target is decomposed into $m$
intervals with the length $\ell_1$, \dots, $\ell_m$. For each interval, we
may assume that line $2$ of $\aw$ runs for $f_h$ times, and that
line $14$ is executed for $g_h$ times. Since $\aw$ will stop (line
$4$) if $T_{start}=T$ ( line $3$), the remainder of $\aw$,
i.e.~lines $7$~to $41$ run for $(f_h - 1)$ times, each such execution
having complexity $\bigo(\ell_h)$. Therefore we arrive at the
complexity
\[
  \sum_{h=1}^m\bigl((f_h+g_h)\crosscomp(\ell_h)+(f_h-1)\bigo(\ell_h)\bigr)\,,
\]
where $\crosscomp(\ell)$ denotes the complexity of the $\foldpk$. The
multiplicities $f_h$ and $g_h$ depend on various factors, such as $start$, the
random order of the elements of $U_1$,$U_2$ (see Algorithm~\ref{A:sls}) and
the probability $p$.
According to \cite{Fenix:08} the complexity of $\crosscomp(\ell_h)$ is
$\bigo(e^{1.146\,\ell_h})$ and accordingly the complexity of ${\tt Inv}$ is given
by
\[
  \sum_{h=1}^m\bigl((f_h+g_h)\bigo(e^{1.146\,\ell_h})\bigr)\,.
\]
In Figure~\ref{F:time} we present the average inverse folding time
of several natural RNA structures taken from the PKdatabase \cite{W:Database}.
These averages are computed via generating $200$ sequences of the target's
neutral networks. In addition we present in
Table~\ref{T:time_example} the total time for $100$ executions of {\tt Inv}
for an additional set of RNA pseudoknot structures.
%%%
%%%%%%%%%%%%%%%%%%%%%%%%%%%%%%%%%%%%%%%%%%%%%%%%%%%%%%%%%%%%%%%%%%%%%%%%%%%%%%%%
%%%
\begin{figure}
  \centering
  \includegraphics[width=\columnwidth]{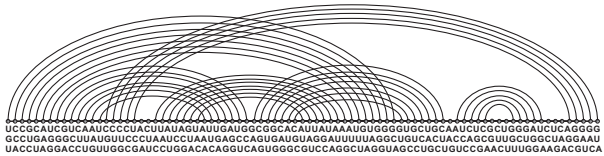}
  \caption{\small A nonplanar $3$-noncrossing RNA structure together with
  three sequences realizing them as mfe-structures.}
  \label{F:pk}
\end{figure}
%%%
%%%%%%%%%%%%%%%%%%%%%%%%%%%%%%%%%%%%%%%%%%%%%%%%%%%%%%%%%%%%%%%%%%%%%%%%%%%%%%%%
%%%
%%%
%%%%%%%%%%%%%%%%%%%%%%%%%%%%%%%%%%%%%%%%%%%%%%%%%%%%%%%%%%%%%%%%%%%%%%%%%%%%%%%%
%%%

\begin{figure}
  \centering
  \includegraphics[width=\columnwidth]{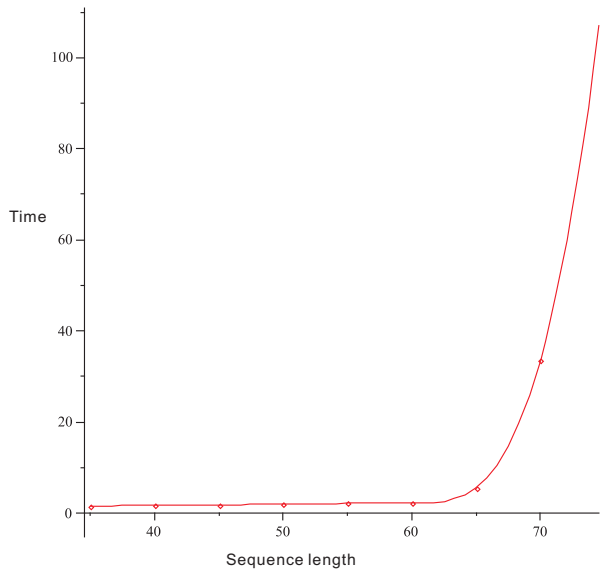}
  \caption{\small
   Approximation using $2$ cubic spines fitting of mean inverse
  folding time (seconds) over sequence length.
  For $n=35,\dots,75$
  we choose a natural pseudoknot structure from the PKdatabase and display
  the average inverse folding time based on sampling $200$ sequences of the
  neutral network of the respective target.}
  \label{F:time}
\end{figure}
%%%
%%%%%%%%%%%%%%%%%%%%%%%%%%%%%%%%%%%%%%%%%%%%%%%%%%%%%%%%%%%%%%%%%%%%%%%%%%%%%%%%
%%%

%%%
%%%%%%%%%%%%%%%%%%%%%%%%%%%%%%%%%%%%%%%%%%%%%%%%%%%%%%%%%%%%%%%%%%%%%%%%%%%%%%%%
%%%
\section{Competing interests}
%%%
%%%%%%%%%%%%%%%%%%%%%%%%%%%%%%%%%%%%%%%%%%%%%%%%%%%%%%%%%%%%%%%%%%%%%%%%%%%%%%%%
%%%

The authors declare that they have no competing interests.

%%%
%%%%%%%%%%%%%%%%%%%%%%%%%%%%%%%%%%%%%%%%%%%%%%%%%%%%%%%%%%%%%%%%%%%%%%%%%%%%%%%%
%%%
\section{Authors contributions}
%%%
%%%%%%%%%%%%%%%%%%%%%%%%%%%%%%%%%%%%%%%%%%%%%%%%%%%%%%%%%%%%%%%%%%%%%%%%%%%%%%%%
%%%

All authors contributed equally to this paper.

%%%
%%%%%%%%%%%%%%%%%%%%%%%%%%%%%%%%%%%%%%%%%%%%%%%%%%%%%%%%%%%%%%%%%%%%%%%%%%%%%%%%
%%%
\section{Acknowledgments}
%%%
%%%%%%%%%%%%%%%%%%%%%%%%%%%%%%%%%%%%%%%%%%%%%%%%%%%%%%%%%%%%%%%%%%%%%%%%%%%%%%%%
%%%
\ifthenelse{\boolean{publ}}{\small}{}

We are grateful to Fenix W.D.~Huang for discussions.
Special thanks belongs to the two anonymous referee's whose thoughtful comments
have greatly helped in deriving an improved version of the paper.
This work was supported by
the $973$ Project, the PCSIRT of the Ministry of Education, the Ministry of
Science and Technology, and the National Science Foundation of China.

%%%
%%%%%%%%%%%%%%%%%%%%%%%%%%%%%%%%%%%%%%%%%%%%%%%%%%%%%%%%%%%%%%%%%%%%%%%%%%%%%%%%
%%%
{\ifthenelse{\boolean{publ}}{\footnotesize}{\small}%
\raggedright%
\bibliographystyle{bmc_article}  % Style BST file
\bibliography{inv}}              % Bibliography file (usually '*.bib' )
%%%
%%%%%%%%%%%%%%%%%%%%%%%%%%%%%%%%%%%%%%%%%%%%%%%%%%%%%%%%%%%%%%%%%%%%%%%%%%%%%%%%
%%%
\section{Tables}
%%%
%%%%%%%%%%%%%%%%%%%%%%%%%%%%%%%%%%%%%%%%%%%%%%%%%%%%%%%%%%%%%%%%%%%%%%%%%%%%%%%%
%%%

\subsection{Table \ref{T:time_example}}

\begin{table}
\begin{tabular}{ l c c r r}
\hline	
  RNA structure & length & trials & total time & success rate	\\ \hline
  TPK-70.28\cite{W:TPK} & 40 & 100 & 4m 57.81s &100\%\\ \hline
  Ec\_PK2\cite{W:Ec} & 59 & 100 & 5m 33.28s &100\% \\
  PMWaV-2\cite{W:PMW}   & 62 & 100 & 1m 7.12s & 100\% \\ \hline
 tRNA & 76 & 100 & 5m 2.49s &100\%\\ \hline
\end{tabular}
\caption{\small Inverse folding times for $100$ executions of {\tt Inv}
for various RNA pseudoknot structures. In all cases all trials generated
successfully sequences of the respective neutral networks.}
\label{T:time_example}
\end{table}
\ifthenelse{\boolean{publ}}{\end{multicols}}{}

\onecolumn % remove when submit

\end{bmcformat}
\end{document}